\documentclass{article}
\usepackage{amsmath}
\font\tengoth=eufm10
\font\sevengoth=eufm7
\font\fivegoth=eufm5
\newfam\gothfam
\textfont\gothfam=\tengoth \scriptfont\gothfam=\sevengoth
\scriptscriptfont\gothfam=\fivegoth

\def\blacksquare{\hbox to .60em{\vrule width .60em height .60em}}

  \font\bb=msbm10 
\catcode`\@=12  

\def\é{\'e}
\def\{\`e}
\def\?{\`a}
\def\{\`u}
\def\{\c c}
\def\hb {\hfil \break}
\def\n {\vskip 0.2cm \noindent }
\def\scirc{\,{\raise 0.8pt\hbox{$\scriptstyle\circ$}}\,}
\def\ins{\,{\raise 0.2cm \hbox{ $\scriptstyle \circ$}}\,}

\def  \é{\'e}
\def\è{\`e}
\def\à{\`a}
\def\ù{\`u}
\def\ç{\c c$\!\!\!$}

\date{}
 
\begin{document}
  
 \centerline {\bf  Calcul explicite  de la courbure des  tissus calibr\és ordinaires :}
 
 \bigskip
 
 \centerline {\bf   J.P. Dufour et D. Lehmann}  
  \bigskip

 \bigskip
 
\section{Introduction } 

\indent  Notons

 $c(n,h):= \frac{(n-1+h)!}{(n-1)!h!}$ la dimension de l'espace vectoriel  des polyn\^omes homog\ènes de degr\é $h$ \à coefficients complexes, en $n$ variables,
 
 $d$ et $n$ deux  entiers ($n\geq 2, \ d>n$), 
 
 $k_0$ l'entier tel que $c(n,k_0)\leq d<c(n,k_0+1)$.
 
 Nous avons d\éfini  dans [CL] les tissus holomorphes de codimension un  dits \emph{ordinaires} (cette terminologie fait allusion au fait que cette propri\ét\é est  g\én\ériquement v\érifi\ée localement). Nous  avons d\émontr\é 
 
 - que le rang d'un tel $d$-tissu   sur une vari\ét\é complexe $M$ de dimension $n$ \était major\é par l'entier $\pi'(n,d)=\sum_{h=1}^{k_0}   \bigl(   d-c(n,h) \bigr)$ (strictement inf\érieur, pour $n\geq 3$, au nombre de Castelnuovo $\pi(n,d)$ qui est  le genre arithm\étique maximum des courbes alg\ébriques irr\éductibles  de degr\é $d$ dans l'espace projectif complexe $\hbox{\bb P}_n$ de dimension $n$),
 
 - que, si  $d$ \était \égal \à $c(n,k_0) $ (nous dirons alors que le tissu est \emph{calibr\é}),     les relations ab\éliennes du tissus au dessus du compl\émentaire $M'=M\setminus S$ d'un certain sous-ensemble analytique $S$ de dimension au plus $n-1$, s'identifiaient  aux sections holomorphes $s$  d'un certain fibr\é holomorphe ${\cal E}\to M'$ de rang $\pi'(n,d)$,  dont la  d\ériv\ée covariante $\nabla s$, relativement \à une certaine connexion holomorphe naturelle $\nabla$ sur ${\cal E}$, \était nulle.
 
\n [Dans le cas $n=2$ des tissus planaires, $S$ est vide, tous les tissus sont \à la fois  calibr\és $\bigl(d=c(2,d-1)\bigr)$ et ordinaires (non-nullit\é des   d\éterminants de Vandermonde). En outre,   $\pi'(2,d)=\pi(2,d)$. La connexion dont il est question a alors \ét\é d\éfinie    par A. H\énaut ([H1]) quel que soit $d$. Sa courbure g\én\éralise celle    donn\ée  par Blaschke-Dubourdieu ([B]) pour $d=3$, et cette g\én\éralisation \à tout $d$ pour $n=2$  avait   \ét\é esquiss\ée  par Pantazzi ([P]), \évidemment dans un autre langage que celui des connexions, qui n'\étaient pas d\éfinies \à cette \époque.]

 La nullit\é de la courbure de cette   connexion  permet donc de montrer que le rang du tissu est     maximum, c'est-\à-dire \égal au rang $\pi'(n,d)$ du fibr\é ${\cal E}$, sans qu'il soit besoin d'exhiber les relations ab\éliennes et de montrer leur ind\épendance lin\éaire. Malheureusement, cette courbure est difficilement calculable ``\à la main", en dehors de cas tr\ès \él\émentaires (m\^eme le cas du 5-tissu de Bol en dimension deux  n\écessite en pratique  l'utilisation d'un ordinateur). 
 
Nous nous \étions    content\és, dans [CL],  de montrer l'existence de cette connexion, en esquissant une m\éthode de calcul, que nous avons essay\é de mettre en oeuvre sur ordinateur avec L. Flaminio et Y. Hantout, mais qui  n'a pas abouti avec eux. Depuis, nous avons r\éussi \à la faire fonctionner, mais   nous proposons dans cet article une autre m\éthode, plus simple, n\écessitant des temps de calcul beaucoup plus courts, avec programmation\footnote{Un programme existerait d\éj\à, r\édig\é par O. Ripoll, mais uniquement dans le cas   $n=2$ des tissus planaires, et  qui n'a pas \ét\é publi\é \à notre connaissance.} explicite des calculs correspondants sur MAPLE.   
Nous expliquerons plus loin  la diff\érence entre ces deux m\éthodes.

\pagebreak

Avec notre programme\footnote{Nous ne sommes pas experts en MAPLE, et notre programme n'est peut-\^etre pas le plus astucieux qu'on puisse imaginer ; mais il a au moins le m\érite d'exister et  de fonctionner.}, nous avons pu 
non seulement  v\érifier la nullit\é de la courbure de certains tissus dont on savait d\éj\à qu'ils \étaient   de rang maximum $\pi'(n,d)$, mais aussi en d\écouvrir d'autres.

Parmi les premiers figurent en particulier le  tissu de Bol   et les $c(n,4)$-tissus  de Pereira-Pirio ($k_0=4$)  qui le g\én\éralisent, not\és ${\cal W}(A_{0,n+3})$ dans [Pe]. Ce sont les tissus dont les int\égrales premi\ères sont les birapports de toutes les familles de 4 nombres (pris dans un ordre arbitraire) parmi les $(n+3)$ nombres   $(x_1,\cdots,x_n,0,1,\infty)$ suppos\és distincts dans la droite projective complexe. Il a \ét\é d\émontr\é que ces tissus \étaient ordinaires et de rang maximum, par Bol pour $n=2$, Pereira-Pirio pour $n=3$, Pereira pour $n$ quelconque ([P]).  Nous avons r\édig\é une feuille de calcul dans laquelle il suffit de pr\éciser $n$, et  avons effectivement fait tourner le programme pour $n=2,3,4$ et 5.

Nous avons aussi v\érifi\é la nullit\é de la courbure du 9-tissu planaire exceptionnel de G. Robert.

Pour r\épondre n\égativement \à une conjecture de Chern-Griffiths selon laquelle, parmi les relations ab\éliennes des  tissus exceptionnels, il y en avait  n\écessairement  qui faisaient intervenir des polylogarithmes, Pirio  a exhib\é ([Pi]) 6 (resp. 10) relations ab\éliennes purement alg\ébriques et lin\éairement ind\épendantes des 5-tissus planaires    $(x,y, x+y, x-y,x^2+y^2 \hbox{ ou } xy)$ $\bigl($resp. du 6-tissu planaire $(x,y,x+y,x-y,x^2+y^2,xy\bigr)$. Nous avons \évidemment v\érifi\é que la courbure de ces tissus \était nulle, ainsi d'ailleurs que celle du  7-tissu $\bigl(x,y,x+y,x-y,x^2+y^2,x^2-y^2,xy\bigr)$ de rang 15. Mais, l'ayant  constat\ée par ordinateur   pour $n=3$ ou 4,   nous   avons pu  en d\éduire   pour tout $n$ la nullit\é de la courbure des $c(n,4)$ tissus 
${\cal W}B_n$  du m\^eme type d\éfinis par  les fonctions 
 $$x_i,\ x_i+x_j,\ x_j-x_i,\ \ x_i.x_j,\ x_i+x_j+x_k, \ x_i^2+x_j^2+x_k^2,\  x_i.x_j.x_k $$    et 
 l'une des familles de   fonctions $$x_i+x_j+x_k+x_m\hbox{ ou }x_i^2+x_j^2+x_k^2+x_m^2\hbox{ ou }x_i.x_j.x_k.x_m .$$(voir la d\émonstration dans la section 7).

 

Le programme est en fait  tr\ès sensible \à l'ordre des int\égrales premi\ères $u_i$. Cela tient \à ce  qu'on lui impose les variables principales dans le calcul d'une trivialisation locale de $\cal E$.
Par exemple, pour le 15-tissu ${\cal W}(A_{0,6})$, on a besoin  d'inverser une sous-matrice carr\ée $YYY$ de taille $19\times 19$ d'une matrice $MM$ de taille $19\times 45$ ; la simple  transposition des birapports $(x,y,z,\infty)$ et  $(x ,y ,0,1)$  conduit alors \à une matrice $YYY$  de rang 18, tandis qu'une   autre permutation particuli\ère des $u_i$ conduit  \à une matrice $YYY$ qui est bien de rang 19, mais avec un temps de  calcul final de 1'46'' au lieu de 28''. 

On pourrait vouloir   rem\édier \à ce  d\éfaut en laissant MAPLE choisir lui-m\^eme les variables principales. Mais on s'aper\ç\  oit alors qu'il fait ce choix de fa\ç on plus ou moins al\éatoire, variant d'une fois \à l'autre,  et que  l'on ne contr\^ole plus rien du tout. En outre, le r\ésultat n'est alors  lisible qu'en cas de courbure nulle, puisque  la trivialisation par rapport \à laquelle est    calcul\ée la  forme de courbure varie ; or cette courbure   est  toujours un invariant int\éressant du tissu : en particulier,   l'existence d'un sous-fibr\é de $\cal E$, invariant par la connexion, et sur lequel la courbure est nulle, permet de borner inf\érieurement le rang du tissu par le rang de ce sous-fibr\é (on observe imm\édiatement par exemple que  le 8-tissu planaire $\bigl(x,y,x+y,x-y,x^2+y^2,x^2-y^2,xy,x^4+y^4\bigr)$  n'est pas de rang maximum 21, mais qu'il est au moins de rang 19). 

  \section{Rappel de   d\éfinitions   des relations ab\éliennes}

  Un $d$-tissu sur $M$ n'est en fait qu'un feuilletage $\tilde{\cal F}$ sur l'espace total   d'un certain rev\^etement   de $M$     
 \à $d$-feuillets.
  Sur un ouvert   $U$ de $M$ au dessus duquel ce rev\^etement est trivial,
  il revient   au m\^eme   
  de se donner 
  $d$ feuilletages $ {\cal F}_i$ sur $M$ : on dit alors que le tissu est \emph{totalement d\écomposable} au dessus de $U$.

  Le calcul de la courbure,  \étant local, permet de se restreindre \à un tel ouvert 
 et nous  nous contenterons, pour simplifier l'expos\é, de rappeler les d\éfinitions dans le cas d'un tissu totalement  d\écomposable.    Mais le fibr\é $\cal E$ dans le cas des tissus ordinaires, et la connexion tautologique dont on le munit dans le cas des tissus qui sont en plus calibr\és, sont en fait d\éfinis gloçbalement sur tout l'ouvert $M'=M\setminus S$.

 On suppose donc le tissu  d\éfini   par la donn\ée de $d$ feuilletages   holomorphes ${\cal F}_i$ de codimension 1 sur la vari\ét\é complexe $M$ de dimension $n$, ($d>n$),     en position g\én\érale au moins faible\footnote{ On dit que le tissu est en position g\én\érale faible (resp. forte)  si, en tout point $m$ de la partie r\éguli\ère $M_0$ du tissu, il existe au moins $n$ des feuilletages parmi les $d$ dont les espaces tangents en $m$ sont en position g\én\érale (resp. si toute famille de $n$ feuilletages  parmi les $d$ a cette propri\ét\é).}. Quitte \à remplacer $M$ par un ouvert, on supposera d\ésormais   le tissu r\égulier  sur tout $M$.

 Une   \emph{relation ab\élienne}  sur un ouvert $U$ (suppos\é connexe et simplement connexe) de $M$ est alors la donn\ée d'une famille $(F_i)_i$ de fonctions holomorphes sur $U$, $1\leq i\leq d$, telles que
 
 - pour tout $i$, $F_i$ est une int\égrale premi\ère de ${\cal F}_i$ (avec \éventuellement des singularit\és)), 
 
 - la somme $\sum_{i=1}^d F_i$ est une fonction constante sur $U$. 
 
 \n Ces int\égrales premi\ères n'\étant d\éfinies qu'\à une constante additive pr\ès, cela revient en fait \à ne d\éfinir que leur diff\érentielle $\omega_i=dF_i$, de sorte que l'on peut encore d\éfinir une relation ab\élienne comme 
  une  famille  $(\omega_i)_i$, $1\leq i\leq d$,  de 1-formes holomorphes $\omega_i$ sur $U$ (admettant \éventuellement des singularit\és), qui sont   
 
  $(i)$ toutes       ferm\ées  (donc localement  exactes)  : $d \omega_i=0$,

   $(ii)$  qui v\érifient  $T{\cal F}_i\subset Ker \  \omega_i$ quel que soit $i$ \  ($T{\cal F}_i=Ker \  \omega_i$ en tout point o\ù $\omega_i$  n'est pas nulle), 
   
   $(iii)$ telle que $\sum_{i=1}^d  \omega_i=0$.

 Dans le but d'introduire $({\cal E},\nabla)$, nous allons en donner une d\éfinition \équivalente en termes d'op\érateurs diff\érentiels.
  Notons $T{\cal F}_i$ le sous-fibr\é de $TM$ form\é des vecteurs tangents \à ${\cal F}_i$, et $A_i$ le sous-fibr\é de $T^*M$ form\é des formes lin\éaires nulles sur $T{\cal F}_i$ (dual du fibr\é quotient $TM/T{\cal F}_i$).
  
  \n Soit   $Tr:\oplus_{i=1}^d A_i\to  T^*M$  l'homomorphisme de fibr\és vectoriels (appel\é \emph{Trace}), d\éfini par   $$ Tr\bigl(( \omega_i)_i\bigr)= \sum_{i=1}^d  \omega_i .$$ L'hypoth\èse de position g\én\érale au moins faible signifie  qu'il est de rang maximum $n$ : son noyau $A=Ker \ Tr$ est donc un fibr\é vectoriel de rang $d-n$.  
  
 On d\éfinit   un op\érateur diff\érentiel lin\éaire $D:J^1 A\to B$   d'ordre 1,     o\ù $B=(\wedge ^2T^* M)^{\oplus d}$, en associant, \à toute section $s=( \omega_i)_i$ de $A$, la  famille $(d \omega_i)_i$ des diff\érentielles. Une  \emph{relation ab\élienne}   est alors une  section holomorphe $s$ de $A$ telle que $D(j^1s)=0$.
 
 
  
  
  Avec les notations pr\éc\édentes, posons  $R_1=Ker (D:J^1 A\to B)   $ : c'est l'espace des \emph{relations ab\éliennes formelles \à l'ordre 1}. Plus g\én\éralement, l'espace des \emph{relations ab\éliennes formelles \à l'ordre h} est le noyau du $(h-1)$-i\ème prolongement  $D_h$ de l'op\érateur diff\érentiel $D\ (=D_1)$ :
  $$R_h=Ker (D_h:J^h A\to J^{h-1}B)  .$$ 
  
  Notant $\pi_h:R_h\to R_{h-1}$ la restriction \à $R_h$ de la projection $J^h A\to J^{h-1} A$, nous  montrerons  que les \él\éments de $R_h$ qui se projettent par $\pi_h$ sur un \é\l\ément donn\é $a_{h-1}$ de $R_{h-1}$ sont les solutions d'un syst\ème lin\éaire $\Sigma_h(a_{h-1})$ de $c(n,h+1)$ \équations \à $d$ inconnues, dont la partie homog\ène   ne d\épend pas de $a_{h-1}$. 
 \section{Comparaison des deux m\éthodes}
 
 Supposons chaque feuilletage ${\cal F}_i$  d\éfini sur $U$ par une int\égrale premi\ère $u_i$ sans singularit\é.
 
 \n Dans la premi\ère m\éthode, nous avions    suppos\é les coordonn\ées locales $x =(x_1,\cdots,x_n)$
  choisies de \hb  fa\ç on que $\frac{\partial}{\partial x_n}$ soit transverse \à  tous les feuilletages   ${\cal F}_i$ 
  ($1\leq i\leq d$) ; ceux-ci pouvaient donc  \^etre d\éfinis par les 1-formes holomorphes int\égrables  
  $$\eta_i=dx_n+\sum_{\alpha=1}^{n-1} p_{i \alpha}( z)\ dx_\alpha,\hbox{   o\ù  } p_{i \alpha}=\frac{(u_i)'_\alpha}{(u_i)'_n}.$$
 Les relations ab\éliennes $(\omega_i)$ \étaient  d\éfinies par les familles $(f_i)$ de fonctions holomorphes telles que $\omega_i=f_i\ \eta_i$, et les relations ab\éliennes  formelles \à l'ordre $k$ en un point $m$ \étaient d\éfinies par la valeur des $f_i$ et de  leurs d\ériv\ées partielles  jusqu'\à  l'ordre $k$ en ce point. Mais nous r\éduisions \à $d. (k+1)$ le nombre de ces inconnues en observant que toutes les d\ériv\ées partielles cherch\ées s'exprimaient \à l'aide des seules d\ériv\ées partielles successives par rapport \à $x_n$.
 
 Maintenant, nous ne faisons plus jouer de r\^ole particulier \à la coordonn\ée $x_n$, et cherchons les int\égrales premi\ères $F_i$ sous la forme $ F_i=G_i(u_i)$, o\ù $G_i$ d\ésigne une fonction holomorphe d'une seule variable. Notant $g_i$ la d\ériv\ée de $G_i$, on a maintenant $\omega_i=g_i(u_i)\ du_i$,  et les relations ab\éliennes  formelles \à l'ordre $k$ en un point $m$ sont alors  d\éfinies par la valeur des fonctions $ g_i$ et de  leurs d\ériv\ées  jusqu'\à  l'ordre $k$ au point $u_i(m)$, de sorte que le nombre d'inconnues est encore \égal \à $d .(k+1)$. Mais, outre que  les coordonn\ées locales jouent d\ésormais  toutes le m\^eme r\^ole, les calculs sont beaucoup plus simples et plus rapides. 
 
\section{Localisation et calculs}
\noindent {\bf Notations }

 \n $i$ d\ésigne un indice variant de 1 \à $d$,
 
 \n $\lambda,\mu,...$ des entiers variant de 1 \à $n$.

 \n  $L=(\ell_1,\ell_2,\cdots,\ell_n)$ d\ésigne un multi-indice form\é d'entiers $\ell_\lambda\geq 0$,\hb   $|L|:=\sum_\lambda \ell_\lambda$ s'appelle le \emph{degr\é} de $L$. 
  
  \n Si $L=(\ell_1,\ell_2,\cdots,\ell_n))$, et $L'=(\ell'_1,\ell'_2,\cdots,\ell'_n))$, $L+L':=(\ell_1+\ell'_1,\cdots,\ell_n+\ell'_n)$.

 \n $ 1_\lambda$ d\ésigne le multi-indice obtenu avec 1 \à la place $\lambda$ et 0 ailleurs.

   \n Lorsque $\ell_\lambda\geq 1$, $L-1_\lambda$ d\ésigne le multi-indice obtenu en rempla\ç ant $\ell_\lambda$ par $\ell_\lambda -1$.

  \n Relativement \à des coordonn\ées locales $x =(x_1,\cdots,x_n)$    dans {\bb C}$^n$,   on notera $\partial_\lambda a$ ou $a'_\lambda$ la  d\ériv\ée partielle
  $\frac{\partial a}{ \partial x_\lambda } $ d'une fonction holomorphe $a$ ou d'une matrice \à coefficients holomorphes.  \n Plus g\én\éralement, $\partial_L a$ ou $a'_L$ d\ésigne la d\ériv\ée partielle $\frac{\partial^{|L|}a}{(\partial x_1)^{\ell_1}\cdots(\partial x_n)^{\ell_n}} $ d'ordre $|L|$.


   \subsection{Principes du calcul}
   
  Se donner une int\égrale premi\ère $F_i=G_i(u_i)$ \à constante additive pr\ès  revient maintenant \à se donner la fonction d\ériv\ée $g_i=(G_i)'$.  Chaque fibr\é $A_i$ \étant d\ésormais trivialis\é par $du_i$, on pose  $\omega_i=g_i(u_i)\ du_i$ (une telle forme est n\éc\éssairement ferm\ée). Se donner une 
  relation ab\élienne revient alors \à se donner  une famille   $(g_i)$ de fonctions holomorphes d'une variable  ($1\leq i\leq d$) telles que $\sum_i g_i(u_i)\ du_i\equiv 0$, 
 soit :  $$  \sum_i g_i(u_i)\ (u_i)'_\lambda \equiv 0\hbox {\hskip 1cm  pour    tout $ \lambda$}.$$
 [Les fonctions $u_i$  sont donn\ées ; les inconnues sont les fonctions $g_i$]. 
 
\n Soit $R_h$ l'espace des relations ab\éliennes formelles \à l'ordre $h$. 
Les   d\ériv\ées partielles successives des relations $(iii)$ vont permettre de calculer   $R_h$.

\subsection{Les coefficients  $M^h_{ L}( u )  $}

\n {\bf Lemme 1 :} {\it Pour toute  fonction  holomorphe $u$ de $n$ variables, et toute fonction holomorphe $g$ d'une variable, 

$(i)$ les d\ériv\ées $\Bigl((g\scirc u)\ u'_\lambda\Bigr)'_L$,  sont  des combinaisons lin\éaires
$$\Bigl((g\scirc u)\ u'_\lambda\Bigr)'_L=\sum_{h=0}^{|L|} M_{L+1_\lambda}^h(u)\ .\ (g^{(h)}\scirc u)  $$des   d\ériv\ées successives $g^{(h)}$ de  $g$ $($on a pos\é $g^{(0)}=g)$, dont  les coefficients  $M_{L' }^h(u) =M_{L+1_\lambda}^h(u)$ ne d\épendent que de $u$   et  du multi-indice $L'=L+1_\lambda$, et non de la d\écomposition de celui-ci sous la forme $L+1_\lambda$.

$(ii)$ On les calcule par r\écurrence sur $|L|$ \à l'aide des formules

 $$\begin{matrix}
  M^0_{1_\lambda} (u)  & =&\ u'_\lambda\ ,&&\\
    &&&&&  \\
 
  M^0_{  L+1_\mu}(u)&=& \partial_\mu M^0_{  L}(u)\ , & &\\
 
  &&&&&  \\
 
  M^h_{  L+1_\mu}(u)&=& \partial_\mu M^h_{  L}(u)\  +&\ M^{h-1}_{ L} (u)  \ .\  u'_\mu  &\hbox{ si }  1\leq  h\leq |L|-1\ ,\\

  &&&&&  \\

M^{|L|}_{L+1_\mu}(u)&=& M^{|L|-1 }_{ L} (u)  \ .\  u'_\mu \ .&& & \\

\end{matrix}  $$ 
 \n En particulier, on obtient les formules  $$M^0_{L} (u)  =\partial_L u \hbox{,\hskip.5cm et \hskip.5cm }
  M^{|L|-1}_{L}(u)= \prod _{\lambda=1}^n  (u'_\lambda)^{\ell_\lambda}  \hbox{ lorsque }L=(\ell_1,\ell_2,\cdots,\ell_n).$$}

\n  Le lemme  r\ésulte imm\édiatement de ce que, la forme $d\bigl(G(u)\bigr)$ est ferm\ée, $G$ d\ésignant une primitive de $g$.

\n Pour un tissu d\éfini localement par les int\égrales premi\ères $u_i$, on posera :  $$M^h_{i,  L}=M^h_{ L}(u_i) .$$

  \subsection{Les \équations $E_L$ }
  
  Il r\ésulte des consid\érations pr\éc\édentes qu'un \él\ément de  $R_h$ au dessus d'un point $m\in M$ est repr\ésent\é par ses composantes  $(j^h_{u_i(m)} g_i)_i\in \oplus_i J^h A_i$, et que chaque composante $j^h_{u_i(m)} g_ii$ est enti\èrement d\éfinie par la famille des nombres  $ \Bigl( w_i^{(k)}=(g_i^{(k)}\scirc u_i)(m)\Bigr)_{0\leq k\leq h} $.    Pour  qu'une telle famille appartienne \à  $R_h$, il faut et il suffit que soient v\érifi\ées toutes les \équations 
  $(E_L)$   pour $1\leq |L|\leq h+1$, d\éfinies de la fa\ç on suivante : si $L=(\ell_1,\cdots, \ell_n)$ avec $|L|\geq 1$, on choisit un indice $\lambda$ tel que $\ell_\lambda\geq 1$ ;  $(E_L)$ d\ésigne alors l'\équation 
 $$(E_L)\hskip 1cm  \sum_{i=1}^d\sum_{h=0}^{|L|-1} M_{iL}^h\ .\ w_i^{(h)}  =0,                        $$
qui ne d\épend pas du choix de l'indice $\lambda$ tel que $\ell_\lambda\geq 1$.
 On en d\éduit le 
 \n {\bf Th\éor\ème 2 :} {\it
 Si un \él\ément $a_{h-1}\in R_{h-1}$ est d\éfini par les nombres $ \Bigl( w_i^{(k)} \Bigr)_{i,k}  $, $(0\leq k\leq h-1)$, les \él\éments de  $  R_{h }$ se projetant sur  $a_{h-1}$ sont  d\éfinis par les $d$ nombres $(w_i^{(h)})_i$ solutions du syst\ème  lin\éaire $\Sigma_h(a_{h-1})$ des  $c(n,h+1)$ \équations 
 $$ (E_L)\hskip 1cm  \sum_i  C_i^L\ .\ w_i^{(h)}=- \sum_{i=1}^d\sum_{k=0}^{h-1} M_{i L}^k\ .\ w_i^{(k)} ,\ \ \ \ (|L|=h +1)$$\à $d$ inconnues $w_i^{(h)}, \ (1\leq i\leq d)$, 
 o\ù $C_i^L=\prod _{\lambda=1}^n  \bigl((u_i)'_\lambda\bigr)^{\ell_\lambda}  \hbox{ lorsque }L=(\ell_1,\ell_2,\cdots,\ell_n)$.}
 
  \n    On d\éfinit, pour tout $1\leq h\leq k_0$, les matrices $P_h=(\!(C_i^L)\!),\ (1\leq i\leq d,\ |L|=h)$,  
   de taille $d\times c(n,h)$, dont les colonnes sont index\ées par les indices $i\in \{1,\cdots,d\}$, et les lignes  par les multi-indices $L$ de degr\é   $|L|=h$, avec  $C_i^L=\prod _{\lambda=1}^n  \bigl((u_i)'_\lambda\bigr)^{\ell_\lambda}  \hbox{ lorsque }L=(\ell_1,\ell_2,\cdots,\ell_n)$.
   
  L'ensemble ${\cal P}(n,h)$ des multi-indices $L$ de degr\é $|L|=h$ est ordonn\é de  
  la fa\ç on  suivante   : la suite $L_1\cdots L_{n-1}$ est l'\écriture d'un entier naturel en base $h+1$ : on  ne conserve que ceux de ces entiers tels que $\sum_{i=1}^{n-1}L_i\leq h$, et  l'on compl\ète la suite $L_1,\cdots,L_{n-1}$ par $L_n:=h-\sum_{i=1}^{n-1}L_i$ : l'ensemble ${\cal P}(n,h)$ est maintenant identifi\é \à un sous-ensemble de {\bb N}, que  l'on  munit de l'ordre induit.

  \n Ceci permet d'\écrire le syst\ème $\Sigma_h(a_{h-1})$ sous la forme matricielle $$<P_{h+1}, w^{(h)}>=S(a_{h-1}),$$
 o\ù $w^{(h)}$ d\ésigne la matrice colonne \à $d$ lignes $(\!(w_i^{(h)})\!)_{_i}$, et $S(a_{h-1})$ la matrice colonne \à $c(n,h+1)$ lignes $(\!(- \sum_{i=1}^d\sum_{k=0}^{h-1} M_{i L}^k\ .\ w_i^{(k)} )\!)_{_{|L|=h+1}}$.
 
 \subsection{Majoration du rang des tissus ordinaires}
 
 Dire que le tissu est \emph{ordinaire} signifie que tous les syst\èmes $\Sigma_h(a_{h-1})$ sont de rang maximum $$inf \bigl(d,c(n,h+1)\bigr)$$lorsque $a_{h-1}$ se projette sur $M$ en dehors d'un sous-ensemble analytique $S$ de dimension au plus $n-1$ (ou \éventuellement vide). 
 
 Notons  $k_0$ l'entier   tel que $$c(n,k_0)\leq d<c(n,k_0+1),$$ Puique  $C_i^{L+1_n}=(u_i)'_nC_i^L$,  il suffit que les syst\èmes $\Sigma_h$ soient de rang maximum $c(n,h+1)$ lorsque $h\leq k_0$,  pour qu'ils soient encore de rang maximum $d$ quel que soit $h>k_0$. 
 Il revient donc encore au m\^eme de dire que le tissu est ordinaire si le syst\ème des equations $E_L$ pour $|L|\leq k_0-1$   est de rang maximum $c(n+1,k_0-1)-1$ et   $P_{k_0}$ de rang $d$.

On en d\éduit le 
  
\n {\bf Th\éor\ème 3 :}  
  {\it Si le tissu est ordinaire,  les solutions de $\Sigma_h(a_{h-1})$  forment   un espace affine de dimension $d-c(n,h+1)$ pour $h\leq k_0 -1$, et la restriction de $R_h$ \à $M'$   est alors un fibr\é vectoriel de rang $\sum_{k=1}^{h+1}\bigl(d-c(n,k)\bigr)$ au dessus de $M'$. }

  Pour $h\geq  k_0$, les syst\èmes  $\Sigma_h(a_{h-1})$ ont 0 ou 1 solution, de sorte qu'un  jet infini de relation ab\élienne  formelle en un point $m\in M'$ $($a fortiori  un germe puisqu'on est dans un contexte analytique$)$ est enti\èrement d\étermin\é par sa projection sur $R_{k_0-1}$. On en d\éduit le 
  
  \n {\bf Th\éor\ème 4 }([CL]) : 
   {\it Le rang d'un tissu ordinaire $($c'est-\à-dire le maximum  de la dimension\footnote{Si le tissu est en position g\én\érale forte, A. H\énaut ([H2]) a d\émontr\é que cette dimension ne d\épendait  pas de $m$.} de l'espace des germes de relation ab\élienne en un point $m\in M)$ est  au plus \égal au rang  $$\pi'(n,d):=\sum_{k=1}^{k_0}\bigl(d-c(n,k)\bigr)$$ du fibr\é  $R_{k_0-1}|_{_{M'}}$ $($restriction de $R_{k_0-1}$
  \à $M')$. }

\n En effet, en un point $m \in M'$, toutes les matrices $P_h$ sont de rang maximum 
$c(n,h+1)$ pour $h\leq k_0-1$ et $d$ pour $h\geq k_0$. 
L'espace affine des solutions de chaque   syst\ème 
$\Sigma_h(a_{h-1})$ a donc la dimension $d-c(n,h+1)$ pour $h\leq k_0-1$. Il est de dimension 0 ou est vide pour $h\geq k_0$. On en d\éduit 
que  le rang du tissu est major\é par $\pi'(n,d)$ en les points de $M'$ ; un raisonnement \él\émentaire de semi-continuit\é prouve qu'il est aussi a fortiori major\é par  $\pi'(n,d)$ en les points de $S$. 

\subsection{Connexion tautologique des tissus calibr\és ordinaires}

Si le tissu est calibr\é  $\bigl(d=c(n,k_0)\bigr)$, et    ordinaire, on suppose  d\ésormais (quitte à se restreindre \à un ouvert partout dense) que tous les systèmes $\Sigma_h$ sont de rang maximum pour $h\leq k_0-1$.

Posons ${\cal E}:=R_{k_0-2}$ : c'est   un fibr\é vectoriel holomorphe de rang $\pi'(n,d)=\sum_{k=1}^{k_0-1}\bigl(d-c(n,k)\bigr)$, la projection $R_{k_0-1}\to {\cal E}$ \étant maintenant un \emph{isomorphisme} de fibr\és vectoriels. 
Notons  $u:{\cal E}\buildrel\cong\over \rightarrow R_{k_0-1}$ l'isomorphisme inverse.
  Le fibr\é      
  $$R_{k_0-1} :=J^1 {\cal E}\cap (J^{k_0-1}A)$$ est l'intersection des fibr\és $J^1 {\cal E}$ et $ J^{k_0-1}A$ dans $J^1 ( J^{k_0-2}A)$. Notant 
  $\iota :R_{k_0-1}\subset J^1 {\cal E}$ l'inclusion naturelle, l'application compos\ée $\iota \scirc u: {\cal E}\to J^1 {\cal E}$ est une scission holomorphe  de la suite exacte
  $$0\to T^*M'\otimes {\cal E} \to J^1 {\cal E}\buildrel{\buildrel{\iota \scirc u}\over\longleftarrow}\over{\longrightarrow}   {\cal E}\to 0  $$
   et d\éfinit par cons\équent  une connexion holomorphe  $\nabla $    sur ${\cal E} $, que nous appellerons   \emph{la connexion tautologique},   dont la d\érivation covariante associ\ée est donn\ée par la formule
  $$\nabla s= j^1 s -\bigl(\iota \scirc u \bigr)(s). $$
  
  Puisque  $\iota \scirc u$ se factorise   \à travers $R_{k_0-1}$, il est \équivalent de dire, pour une section  $\sigma $ de $A$, que $j^{k_0-1}  \sigma$ est une section de $R_{k_0-1}$ ou que $\nabla (j^{k_0-2}  \sigma) $ s'annule :  les relations ab\éliennes au dessus de $M'$ sont donc  les sections holomorphes  $\sigma $ de $A$ telles que  $\nabla (j^{k_0-2}  \sigma)=0 $. 
  
Dire que  la connexion tautologique est  sans courbure \équivaut alors \à dire que  le tissu est   de rang maximum $\pi'(n,d)$ (le rang de ${\cal E}$) au voisinage de tout point de $M$. On a ainsi d\émontr\é :

\n {\bf Th\éor\ème 5 }([CL]) :  {\it 
 
 \n $(i)$ Pour les $c(n,k_0)$-tissus ordinaires, les relations ab\éliennes  s'identifient, par l'application $\sigma\to j^{k_0-2}\sigma$, aux sections holomorphes $s$ de ${\cal E}$   dont la d\ériv\ée covariante  $\nabla s $ par rapport \à la connexion tautologique est nulle.
 
 \n $(ii)$ Le tissu   est alors  de rang maximum $\pi'(n,d) $ ssi la courbure de la connexion tautologique est nulle. 
 }

 \n {\bf Expression explicite de la d\érivation covariante :}
 
   Une section   de  ${\cal E}$ est d\éfinie par la donn\ée des fonctions d'une variable $  s_i ^{(h)} $  pour   $0\leq h\leq k_0-2$, $i=0\cdots,d$.
   Notons  alors  :
   
  $\nabla _\lambda$   la d\ériv\ée covariante par rapport \à $\frac{\partial}{\partial x_\lambda}$\ ,
    
    et $\bigl(U_i(s)\bigr)_{i=1\cdots,d}$ la solution du syst\ème cram\érien $\Sigma_{k_0-1}(s)$.
   La d\éfinition de $\nabla $ s'exprime  alors par les formules :
 $$\begin{matrix}
 (\nabla _\lambda s)_i^{(h)}\scirc u_i&=& \frac{\partial }{\partial z_\lambda}\ ( s_i ^{(h)}\scirc u_i)&-&u'_\lambda\ . \ (s_i^{(h+1)}\scirc u_i)  &  \hbox {\hskip 1cm  pour } h\leq k_0-3,  \cr  &&&&&  \cr
 \hbox{et }(\nabla _\lambda  s)_i^{(k_0-2)}\scirc u_i&=&  \frac{\partial }{\partial z_\lambda}\ ( s_i ^{(k_0-2)}\scirc u_i)&-& u'_\lambda\ . \  U_i(s) .\cr 
 \end{matrix}$$
 
 \n {\bf Remarque :}  Les formules ci-dessus montrent que, localement, la connexion sur $\cal E$ est la restriction d'une connexion d\éfinie sur un fibr\é trivial de rang $(k_0-1)d$ ;  mais celui-ci n'a aucune signification intrins\èque, alors que   $\cal E$  et la connexion tautologique ont une signification intrins\èque globale, ind\épendante des coordonn\ées locales et du  choix des int\égrales premi\ères  locales $u_i$,  le tissu n'ayant m\^eme pas \à  \^etre totalement d\écomposable globalement.
 
  \n {\bf Trivialisation de $\cal E$ et forme de connexion :}
  
  Notons $\pi_h(s)$ la projection sur $R_h$ d'une section $s$ de $\cal E$  pour $h\leq k_0-2$. La famille  $s^{(h)} = (s_i^{(h)})_{_i}$  appartient \à l'ensemble des  solutions de $\Sigma\bigl(\pi_{h-1}(s)\bigr)$ : il suffit donc de d\éfinir, pour tout $ h\leq k_0-2$, un sous-ensemble $I_h$ de $d-c(n,h+1)$ entiers  $i\in \{1,\cdots,d\}$ pour lesquels la sous-matrice correspondante de $P_{h+1}$ (c'est-\à-dire celle dont les colonnes sont index\ées par les  indices $i\notin I_h)$    a un d\éterminant non-nul\footnote{Attention : Il se peut qu'il faille encore restreindre l'ouvert de $M'$ au dessus duquel on se place, car il n'est peut-\^etre pas possible d'utiliser le m\^eme ensemble $I_h$ en tous les points de $M'$.}, pour 
     en d\éduire une trivialisation locale de ${\cal E}$   donn\ée par les sections 
  $S_{i,h}$,  ($i\in I_h$, $0\leq h\leq k_0-2$), ainsi d\éfinies ; dans notre  programmation, nous avons choisi pour $I_h$ les $c(n,h+1)$ derniers indices. C'est  bien-entendu  cel\à qui rend le programme sensible \à l'ordre des $u_i$. 
  
 \n {\it Pour $i\in I_h$, $h\leq k_0-2$,   
 $$  (S_{i,h})_j ^{(k)}\equiv 0 \hbox{ si }(j,k)\neq (i,h), \hskip 1cm (S_{i,h})_i ^{(h)}\equiv 1\ .$$ }
On peut  en particulier calculer     $\nabla (S_{i,h})=\sum_\lambda\bigl(\nabla_\lambda (S_{i,h})\ dz_\lambda\bigr)$, ($h\leq k_0-2, \ i\in I_h$), d'o\ù la forme de connexion $\omega$ relative \à cette  trivialisation de ${\cal E}$. 


 \n {\bf   La courbure   : }

On calcule   alors la courbure 

- soit par  les crochets 
$$K_{ \lambda \mu}=[\nabla_\lambda,\nabla_\mu],$$ \hskip 1cm{  c'est-\à-dire :}
 $$    <K, s>= \sum_{\lambda<\mu}\Bigl(\nabla_\lambda\bigl(\nabla_\mu s\bigr)-\nabla_\mu\bigl(\nabla_\lambda s\bigr)\Bigr)\ dz_\lambda\wedge dz_\mu\ ,$$
 
 - soit par la forme de courbure  $\Omega=d\omega+\frac{1}{2}\omega\wedge \omega$ relative \à la trivialisation pr\éc\édente.

   \section{Programmation sur MAPLE 8}
 Exemple des tissus ${\cal W}(A_{0,n+3})$ de Pereira-Pirio (Bol pour $n=2$),
avec une d\éformation en G que l'on peut supprimer pour n grand afin de ne pas allonger exag\ér\ément les temps de calcul.

$>$ restart;

\n {\bf Entr\ée des param\ètres de base :}

$>$ n:= $\cdots$ ; \ k0:=4 \ ;

                               \n {\it Quelques entiers qu'on en d\éduit, qui seront utiles dans la suite $($l'espace $\cal E$ des $(k_0-2)$ jets de relation ab\éliennes sera un fibr\é de rang ro, inclus dans l'espace des vecteurs de dimension alpha; ses \él\éments seront le noyau de la matrice MM \à beta lignes et alpha colonnes$)$.  L'entier ro est aussi \égal au rang maximum $\pi'(n,d)$ du tissu.}

$>$ d:=binomial(n-1+k0,k0);

    \n {\it  $($cette condition exprime que le tissu est "calibr\é"$)$.}
    
$> $ alpha:=(k0-1)*d;\ ro:=k0*d-binomial(k0+n,k0)+1;\ beta:=alpha-ro;

$>$ X:=[seq(x[i],i=1..n)];

 \n {\bf  Entr\ée des d int\égrales premi\ères }({\it comme fonctions  des op$(j,X)$, j=1..n$)$ :}
 
$>$ apply(u,j,X):

$>$  for j to n do u(j,X):=op(j,X) od;

$>$ for j from 2 to n do for i to j-1 do u(n+i+binomial(j-1,2),X):=op(j,X)/op(i,X);

 print(u(n+i+binomial(j-1,2))=(\%)); od od;

$>$ for j from 2 to n do for i to j-1 do 

u(n+binomial(n,2)+i+binomial(j-1,2),X):=(op(j,X)-1+G)/(op(i,X)-1);

 print(u(n+binomial(n,2)+i+binomial(j-1,2))=(\%)); od od;

$>$  for j from 2 to n do for i to j-1 do

 u(n+2*binomial(n,2)+binomial(n,3)+i+binomial(j-1,2),X)
 :=
 
 op(i,X)*(op(j,X)-1)/(op(j,X)*(op(i,X)-1));
 
 print(u(n+2*binomial(n,2)+binomial(n,3)+i+binomial(j-1,2))=\%); od od;

$>$ for k from 3 to n do for j from 2 to k-1 do for i to j-1 do

 u(n+2*binomial(n,2)+i+binomial(j-1,2)+binomial(k-1,3),X):=(op(i,X)-op(k,X))/(op(j,X)-op(k,X));

 print(u(n+2*binomial(n,2)+i+binomial(j-1,2)+binomial(k-1,3))=(\%)); od od od;

                        $>$ for k from 3 to n do 
                        for j from 2 to k-1 do for i to j-1 do
                        
                         u(n+3*binomial(n,2)+binomial(n,3)+i+binomial(j-1,2)+binomial(k-1,3),X)
                         :=
                         
                         op(j,X)*(op(i,X)-op(k,X))/(op(i,X)*(op(j,X)-op(k,X)));
                         
                         print(u(n+3*binomial(n,2)+binomial(n,3)+i+binomial(j-1,2)+binomial(k-1,3))=\%); od od od;

$>$ for k from 3 to n do for j from 2 to k-1 do for i to j-1 do

u(n+3*binomial(n,2)+2*binomial(n,3)+i+binomial(j-1,2)+binomial(k-1,3),X):=

(op(j,X)-1)*(op(i,X)-op(k,X))/((op(i,X)-1)*(op(j,X)-op(k,X)));

print(u(n+3*binomial(n,2)+2*binomial(n,3)+i+binomial(j-1,2)+binomial(k-1,3))=\%); od od od;

$>$ for m from 4 to n do for k from 3 to m-1 do for j from 2 to k-1 do for i to j-1 do

 u(n+3*binomial(n,2)+3*binomial(n,3)+i+binomial(j-1,2)+binomial(k-1,3)+binomial(m-1,4),X):=
 
 (op(j,X)-op(m,X))*(op(i,X)-op(k,X))/((op(i,X)-op(m,X))*(op(j,X)-op(k,X)));
 
 print(u(n+3*binomial(n,2)+3*binomial(n,3)+i+binomial(j-1,2)+binomial(k-1,3)+binomial(m-1,4))=\%); od od od od;

{\bf A partir de maintenant, le programme ne d\épend plus du tissu introduit.}

{\bf Calcul des coefficients M(j,h,L) : }

 $>$ with(LinearAlgebra):
 
 $>$ interface(rtablesize=(k0)*d);
 
 $>$ apply(M,j,h,L,X):

{\it Calcul des premiers coefficients} M(j,0,EE(i))

 $>$ apply(EE,i):
 
 $>$ apply(delta,t,s):
 
 $>$ for t to n do for s to n do if (t=s) then delta(t,s):=1 else delta(t,s):=0 end if od od;
 
 $>$ for i to n do EE(i):=[seq(delta(i,s),s=1..n)] od;
 
 \n ({\it $EE(\lambda)$ est le multi-indice not\é $1_\lambda$ ci-dessus}).
 
 $>$ for j to d do for i to n do M(j,0,EE(i),X):= diff(u(j,X),x[i]) ;
  $>$  od od;
 
{\it G\én\ération et indexation des multi-indices de d\érivation }{\it  d'ordre 0 \à $k_0$:}

 $>$ apply(L,tau):apply(E,r,y):apply(LL,z):for l from 0 to $k0*(k0+1)^(n-2)$ do for k to k0 do E(l,k):=Vector(n) od od:
 
 $>$ tau:=1:
 
 $>$ for k to k0 do for l from 0 to $k*(k+1)^(n-2)$  do                                                                                                   p:=l:                                                                                                                              for s to n-1 do r:=p mod (k+1);E(l,k)[n-s]:=r;p:=(p-r)/(k+1) od :                                                                  SS:=sum('E(l,k)[u]', 'u'=1..n-1):                                                                                                     if SS<(k+1) then E(l,k)[n]:=k-SS:L(tau):=E(l,k):tau:=tau+1 end if od od:
   

 $>$ for t to binomial(n+k0,k0)-1 do LL(t):=[seq(L(t)[i],i=1..n)] od ;
 
 $>$ M(j,-1,L)=0 et M(j,h,L):=0 pour h $>\ |L|-1$  
 
 ({\it o\ù $|L| =sum\ (L[i],i=1..n) $} ):

 $>$ for j to d do for tau to binomial(n+k0,k0)-1 do M(j,-1,LL(tau),X):=0 od od ;
 
 $>$ for j to d do for tau to binomial(n+k0,k0)-1 do SS:=sum('LL(tau)[i]','i'=1..n);for h from SS to k0 do M(j,h,LL(tau),X):=0 od od od;

{\it Calcul des M(j,h,L) par r\écurrence sur $|L|$ :}

 $>$ for j to d do                                                                                                                    for ss from 2 to k0 do                                                                                                           for tau to binomial(n+k0,k0)-1 do                                                                                                    if (sum('LL(tau)[i]','i'=1..n)=ss) then                                                                                          for r to n do                                                                                                                     if (LL(tau)[r]=0) then else for h from 0 to ss-1  do \hb \indent  M(j,h,LL(tau),X):=\hb simplify(diff(M(j,h,LL(tau)-EE(r),X),x[r])+M(j,h-1,LL(tau)-EE(r),X)*diff(u(j,X),x[r])) od ;\hb  r:=r+n                                                     fi od fi od od od;

{\it On v\érifie le r\ésultat en les imprimant $($on note provisoirement N(j,h,L) =M(j,h,L,X)$)$ ; cette \étape, qui peut utiliser du temps de calcul lorsque $n$ et $k_0$ augmentent, peut \^etre supprim\ée.}
 
 $>$ for j to d do for h from 0 to k0 do for tau to binomial(n+k0,k0)-1 do\hb\indent  print(N(j,h,LL(tau))=M(j,h,LL(tau),X)) od od od ;

{\bf Les matrices MM, QQ}\ {\it $($matrices des syst\èmes d'\équations $E_L$ pour $|L|\leq k_0-2$ et $|L|= k_0-1$ } {\bf et PP  }$(${\it matrice not\ée $P_{k_0}$ ci-dessus}$)$ :

\n {\it Le tissu est ordinaire si le rang de MM est beta et si celui de PP est  d ; on ne v\érifie pas directement la premi\ère condition, car on aura besoin  d'expliciter ci-dessous une sous-matrice carr\ée YYY de MM, de rang beta.}

{\it Les lignes sont num\érot\ées par l'indice tau de LL(tau). On num\érote maintenant les colonnes :}

 $>$ hh:=j- $>$floor((j-1)/d):ii:=j- $>$j-d*floor((j-1)/d):
 
 $>$ ff:=(tau,eta)- $>$M(ii(eta),hh(eta),LL(tau),X);
 
 $>$ MM:=simplify(Matrix(beta,alpha,ff));

 $>$ fff:=(tau,eta)- $>$M(ii(eta),hh(eta),LL(tau+  binomial(n+k0-1,k0-1)-1),X):
 
 $>$ QQ:=simplify(Matrix(d,(k0-1)*d,fff));
 
 $>$ ffff:=(tau,eta)- $>$M(eta,k0-1,LL(tau+ 
 binomial(n+k0-1,k0-1)-1),X):
 
 $>$ PP:=Matrix(d,d,ffff);

{\bf Calcul d'une base W(j) , j=1..ro, de l'espace des sections de $\cal E$ := Ker (MM) }

{\it D\éfinition d'une sous-matrice carr\ée YYY de MM :}

 $>$ Y(0):=MM:
 
 $>$ for k from 1 to (k0-1) do Y(k):=DeleteColumn  (Y(k-1),(k0-k-1)*d+\hb binomial(k0-k+n-1,n-1)+1..(k0-k)*d) end do:
 
 $>$ YYY:=Y(k0-1);

 $>$ evalb(Rank(YYY)=beta);
 
 ({\it 
si le rang de YYY est strictement inf\érieur \à beta, r\éessayer en modifiant l'ordre des   $u(i))$}.
  
 $>$ IYYY:=simplify(MatrixInverse(YYY));

 $>$ B(0):=MM:
 
 $>$ for k from 1 to  (k0-1) do 
 \hb B(k):=DeleteColumn(B(k-1),(k0-k-1)*d+1..(k0-k-1)*d+binomial(k0-k+n-1,n-1)) end do;
 
 $>$ B:=simplify(B(k0-1));
 
 $>$ for j from 1 to ro do ColB(j):=Column(B,j) end do;
 
 $>$ apply(a,j,s);
 
 $>$ for j from 1 to ro do for s from 1 to beta do a(j,s):=factor(factor\hb (simplify((simplify(-IYYY.ColB(j)))[s])));print(`a(` ,j,s , `)=` , a(j,s)) end do end do;

{\bf Partition de la suite 1..(ko-1)d en R(h) et S(h) pour h de 0 \à $(k_0-2)d$ ;\hb 
nr(h) := nombre d'\él\éments dans R(h) ; ns(h) := nombre d'\él\éments dans S(h) :}

 $>$ apply(R,h);apply(S,h);
 
 $>$ for h from 0 to k0-2 do R(h):=[seq(i,i=h*d+1..h*d+binomial(h+n,n-1))] end do;
 
 $>$ for h from 0 to k0-2  do nr(h):=binomial(h+n,n-1) end do;
 
 $>$ for h from 0 to k0-2  do S(h):=[seq(i,i=h*d+binomial(h+n,n-1)+1..(h+1)*d)] end do;
 
 $>$ for h from 0 to k0-2  do ns(h):=d-binomial(h+n,n-1) end do;
 
 $>$ VV:=proc(j) global V; V:=Vector((k0-1)*d);
  
 $>$ for i in R(0) do V[i]:=a(j,i) od:
 
 $>$ for h from 1 to (k0-2) do for i in R(h) do V[i]:=a(j,i-sum(ns(kkk),kkk=0..h-1)) od; od;
   
 $>$ for h from 0 to (k0-2) do for i in S(h) do
  
 $>$ if i= (j+sum(nr(kk),kk=0..h)) then  V[i]:=1 else V[i]:=0 fi ;od;od;
 
 $>$ evalm(V): end proc ;
 
 $>$ apply(W,j);
 
 $>$ for j to ro do VV(j):W(j):=V od;

{\bf Matrice  U exprimant les termes de rang $k_0-1 $ en fonction des termes de rang inf\érieur:}

 $>$ evalb(Rank(PP)=d);
 
 ({\it deuxi\ème condition pour que le tissu soit ordinaire})

 $>$ IPP:=MatrixInverse(PP);

 $>$ U:=simplify(-IPP.QQ);
 
 $>$ U:=simplify(U);

{\bf D\éfinition de la connexion tautologique sur $\cal E$ : }

{\it Expression des d\ériv\ées covariantes dans l'espace des vecteurs de dimension alpha :}

 $>$ Nabla:=proc(VV,j) description "calcul du   $Nabla_{x(j)}$ du vecteur VV ; le r\ésultat est le vecteur Vec";global Vec; Vec:= Vector((k0-1)*d);
 
 $>$ for h from 1 to (k0-2) do
 
 $>$ for i   to  d do
 
 $>$ Vec[i+(h-1)*d]:=diff(VV[i+(h-1)*d],x[j])- VV[i+h*d]*diff(u(i,X),x[j])od od; 
 
 $>$ for i   to  d do
  Vec[i+(k0-2)*d]:=simplify(diff(VV[i+(k0-2)*d],x[j]) -(U.VV)[i]*diff(u(i,X),x[j]))  od;  end proc; 
  
{\it Calcul  de la forme de connexion  relative \à la trivialisation $(W_j)$ :} 

 $>$  DerCov:=proc(i); "le r\ésultat est la matrice A";  A:=Matrix(ro,ro); 
 
 $>$ apply(DC,j,i):for j to ro do Nabla(W(j),i);DC(j,i):=Vec od:apply(f,i):f(i):=(s,j)-$>$ DC(j,i)[s]:\hb
 \indent apply(N,i):N(i):=Matrix(alpha,ro,f(i));
 
 $>$ apply(Aa ,i,j):Aa(i,0):=N(i):for k to  (k0-1) do \hb  \indent Aa(i,k):=DeleteRow(simplify(Aa(i,k-1)),(k0-k-1)*d+1..(k0-k-1)*d+binomial(k0-k+n-1,n-1)) end do:\hb \indent A:=simplify(Aa(i,k0-1));end proc;

 $>$ apply(A,i):for i to n do
 A(i):=DerCov(i):print(connexion(i)=DerCov(i)) od; 
 
\n {\it $connexion(i)$ est la composante sur $dx_i$ de la forme de connexion    relative  \à la trivialisation $\bigl(W(j)\bigr)_j$ .}
  
{\bf Calcul de la forme de courbure K :}

 $>$ apply(A,r,s): apply(f,r,s):
  
 $>$ for s to n do
   for r to n do
  f(r,s):=(i,j)- $>$simplify(simplify(diff(A(r)[i,j],x[s]))):
  
 $>$ A(r,s):=Matrix(ro,ro,f(r,s)):od od:
   
 $>$ apply(AA,r,s):for s  to n do
 for r to n do
   AA(r,s):=simplify(simplify(A(r).A(s))): od od:
 
 $>$ apply(ff,r,s):
 
 $>$ apply(K,r,s):
 
 $>$ for s from 2 to n do 
  for r to s-1 do
  
  ff(r,s):=(i,j)-$>$simplify(simplify(A(r,s)[i,j]-A(s,r)[i,j]+AA(s,r)[i,j]-AA(r,s)[i,j]));
  
   K(r,s):=Matrix(ro,ro,ff(r,s)):print(courbure(r,s)=K(r,s)) od od ;

{\it Pour les tissus \à param\ètre G dont la courbure s'annule pour G=0, d\éveloppement limit\é en G \à l'ordre 0 de la courbure et  le programme affiche  les \él\éments en $O(G)$.  $($s'il n'y a pas de param\ètre, le r\ésultat est le m\^eme qu'\à la ligne pr\éc\édente$)$. Attention \à ce que que la lettre $G$ n'aie pas \ét\é utilis\ée comme symbole par ailleurs .  }
 
 $>$ apply(ffo,r,s):
 
 $>$ apply(Ko,r,s):
 
 $>$ for s from 2 to n do 
  for r to s-1 do
  ffo(r,s):=(i,j)-$>$  taylor(K(r,s)[i,j],G,1);
  
 $>$ Ko(r,s):=Matrix(ro,ro,ffo(r,s)):print(courbure(r,s)=Ko(r,s)) od od ;
 
 \n {\it courbure$(r,s)$ est la composante sur $dx_r\wedge dx_s$ de la forme de courbure relative \à la trivialisation $\bigl(W(j)\bigr)_j$}
 
   \section{Exemples tests   : }  
   Pour la satisfaction de l'oeil, afin de nous assurer que l'on n'obtenait   pas syst\ématiquement une courbure nulle \à la suite d'une erreur de programmation, nous avons d\éform\é certains des   tissus de ces exemàples \à l'aide d'un param\ètre scalaire $G$ (le tissu concern\é \étant obtenu pour $G=0$), et nous  avons  mis en \évidence  les coefficients non nuls de la courbure qui sont  petits d'ordre G. 
   
   \n Pour   $n=2,3,4$, on note parfois dans cette section    $ (x,y) $, $ (x,y,z) $ , $ (x,y,z,t) $ les  coordonn\ées locales. 
   
\n  1) $n=2$, $k_0=2$  : d\éformation du 3-tissu hexagonal   avec d\éveloppement limit\é en $G$ de la courbure de Blaschke \à l'ordre que l'on veut (temps de calcul : 07'', rang 1) : 
 $$\begin{matrix}  u_1 =&x& &  u_2 =&y & & u_3  =&x+y+Gf(x,y). \\ 
 
\end{matrix} $$ 
 \n   2) $n=3$, $k_0=2$  : l'un des rares exemples, en dehors du cas $(n=2, k_0=2)$ de la courbure de Blaschke,   que l'on peut  traiter sans ordinateur ; cf [CL]); (11'', rang 3);\hb on introduit une fonction holomorphe arbitraire $y\ -\!\!\!>F(y), \ F\neq  25$ :
 $$\begin{matrix}  u_1 =   z &     u_2 =  x+ y+z &     u_3  =  2x+4y+z &
                   u_4 =  3x+9y+z&   u_5  =  4x+16y+z   &  u_6=5x +F(y)+z .\\

 \end{matrix} $$ 
 
 \n 3) Les tissus ${\cal W}(A_{0,n+3})$ de Pereira-Pirio ($n\geq 2$, $k_0=4$), 
avec une d\éformation en G que l'on peut supprimer \à volont\é, en particulier pour n ou $k_0$ grands, afin de ne pas allonger exag\ér\ément les temps de calcul : c'est l'exemple r\édig\é dans la section pr\éc\édente.

\n $n=2$ (tissu de Bol,  rang 6) ; temps de calcul : 10",\hb
$n=3$ (rang 26) ;  temps de calcul : 30" sans d\éformation, 44" avec,\hb
$n=4$ (rang 71) ;  temps de calcul : 4'12" sans d\éformation, 12'02" avec,\hb
$n=5$ (rang 150) ; temps de calcul :  55' sans d\éformation.

 \n 4) G\én\éralisation ${\cal W}B_n$ \à tout $n$ des 5-tissus planaires $(x,y, x+y, x-y,x^2+y^2 \hbox{ ou } xy)$ de Pirio,  qui sont de rang maximum   ($n\geq 2$, $k_0=4$) : voir la section suivante.
 
  \n   5) $n=2$, $k_0=5$  : d\éformation du 6-tissu
 $(x,y, x+y, x-y,x^2+(1+G) y^2 , xy)$ de Pirio (11'', rang 10) , 
 ou du  7-tissu ($k_0=6$, 19", rang 15)  obtenu en ajoutant encore     $u_7=x^2-(1+G)y^2 $.
 \n [quant au 8-tissu    obtenu en ajoutant encore $u_8=x^4+y^4$   lorsque $G=0$, on observe imm\édiatement   que sa courbure n'est pas nulle, mais que le sous-fibr\é engendr\é par les 19 premiers vecteurs de la trivialisation est pr\éserv\é par la connexion et que la restriction de la courbure \à ce  sous-fibr\é est nulle (24") : il est donc de rang 19 ou 20]. 
  
  \n   6) $n=2$, $k_0=8$ : le 9-tissu exceptionnel de G. Robert (5'02" sans d\éformation,  rang 28) :
 $$x,y,\frac{x}{1+y},\frac{1+x}{y}, \frac{x}{y}, \frac{1+x}{1+y}, \frac{y(1+x)}{x(1+y)}, \frac{(1+x)(1+y)}{xy},\frac{x(1+x)}{y(1+y)}.$$
  
 \section{ Les tissus  ${\cal W}B_n$ : }  
 
 On observe que $c(n,4)$ est \égal \à $n+3\frac{n(n-1)}{2}+3\frac{n(n-1)(n-2)}{6}+\frac{n(n-1)(n-2)(n-3)}{24}$ pour $n\geq 4$, \à $15$ pour $n=3$ et \à $5$ pour $n=2$.

\n {\bf Th\éor\ème 6 : }
 
{\it Les $c(n,4)$-tissus ${\cal W}B_n$ admettant  comme int\égrales premi\ères locales les fonctions

 \pagebreak
 \n $x_i,\hb \ x_i+ x_j\ ,\ x_j-x_i\  ,\ \ x_i.x_j\ ,\ (i<j),\hb 
 \ x_i+x_j+x_k\ , \ x_i^2+x_j^2+x_k^2\ , \  x_i.x_j.x_k \ ,(i<j<k)\hb     
x_i.x_j.x_k.x_m, \ (i<j<k<m)$ \hb 

\noindent sont ordinaires et de rang maximum $\pi'\bigl(n,c(n,4)\bigr)$.}

 \n {\it D\émonstration :}  Pour $n=2$,   le r\ésultat est d\émontr\é par Pirio, qui explicite 6 relations ab\éliennes (alg\ébriques) ind\épendantes. Pour $n=3$ ou 4, on v\érifie avec notre programme que  la courbure est nulle (temps de calcul avec d\éformation  : 24" pour $n=3$,  3'19"  pour $n=4$).
Puisque le 5-tissu planaire ${\cal W}B_2$ de Pirio est de rang maximum 6, le tissu ${\cal W}B_n$ contient   exactement $6\ \frac{n!}{2!(n-2)!}$ relations ab\éliennes ind\épendantes ne faisant intervenir que 2 des $n$ coordonn\ées locales. Puisque le 15-tissu ${\cal W}B_3$ a une courbure nulle, c'est que le nombre de relations ab\éliennes ind\épendantes  faisant intervenir les 3   coordonn\ées locales  est \égal \à $8$ $\bigl(=\pi'(3,15)-6 \ \frac{3!}{2!1!}\bigr)$ ; par cons\équent  le nombre  des relations ab\éliennes ind\épendantes de ${\cal W}B_n$ faisant intervenir exactement  3 des $n$   coordonn\ées locales  est \égal \à $8\ \frac{n!}{3!(n-3)!}$. De m\^eme, puisque  le 35-tissu ${\cal W}B_4$ a une courbure nulle, c'est que le nombre de ses relations ab\éliennes ind\épendantes faisant intervenir les 4   coordonn\ées locales  est \égal \à $3$ $\bigl(=\pi'(4,35)-6 \ \frac{4!}{2!2!}-8  \frac{4!}{3!1!}\bigr)$. Il en r\ésulte que, pour $n\geq 4$, ${\cal W}B_n$ poss\ède au moins
$$6\ \frac{n!}{2!(n-2)!}+ 8\ \frac{n!}{3!(n-3)!} +3\ \frac{n!}{4!(n-4)!} $$
  relations relations ab\éliennes ind\épendantes  (ne faisant intervenir que 2, 3 ou 4 variables). Or ce  nombre est pr\écis\ément \égal \à $\pi'\bigl(n,c(n,4)\bigr)$, comme on le v\érifie ais\ément. Il suffit donc de d\émontrer que le tissu est ordinaire pour en d\éduire qu'il est de rang maximum.

Commen\ç ons  par montrer que la matrice $P_4$ ($=PP$) est inversible.

\n - On choisit un ordre  arbitraire O2 sur les couples (i,j) tels que $i<j$ , un ordre arbitraire O3 sur les triplets $(i,j,k)$ tels que $i<j<k $ et un ordre arbitraire O4 sur les quadruplets $(i,j,k,m)$ tels que $i<j<k<m$.
\n - On ordonne les colonnes (qui correspondent aux fonctions $u(i)$)   en mettant d'abord les int\égrales premi\ères $x_1,...,x_n$ puis celles qui font intervenir deux variables dans l'ordre O2,  puis celles qui font intervenir trois variables dans l'ordre O3, puis les fonctions  $x_i.x_j.x_k.x_m$ dans l'ordre O4.

\n - On r\éordonne maintenant les lignes (qui correspondent aux multi-indices $L$ d'ordre 4).
On met en premier ceux qui ne font intervenir qu'une variable $4_1,...,4_n$ ($j_i$ d\ésignant le multi-indice dont tous les termes sont nuls sauf le i-\ème \égal \à $j$) ; on range ensuite  ceux qui font intervenir deux variables 
$ x_i$ et $x_j$ en rangeant  les triplets $(3_i+1_j, 2_i+2_j,1_i+3_j)$ pour $i<j$ suivant O2 ; 
 on range ensuite ceux qui font intervenir les trois variables $ x_i$, $ x_j$ et $x_k$ en ordonnant les triplets $( 2_i+1_j+1_k, 1_i+2_j+1_k, 1_i+1_j+2_k)$ pour $i<j<k $ suivant O3 ;
  on range enfin les $1_i+1_j+1_k+1_m$ pour $i<j<k<m$ suivant O4.
  
\n Apr\ès ces r\é-ordonnancements $P_4$ devient une  matrice par blocs, dont tous les   blocs diagonaux sont inversibles, et les blocs sous les blocs diagonaux sont nuls. 

Pour $P_3$ (resp. $P_2$, resp. $P_1$),  on agit de m\^eme, mais en ne conservant que la sous-matrice obtenue en \éliminant les colonnes correspondant aux $x_i.x_j$, $x_i^2+x_j^2+x_k^2$, $x_i.x_j.x_k$ et $x_i.x_j.x_k.x_m$ (resp.  en ef-\hb fa\ç ant toutes les colonnes sauf celles qui correspondent aux $x_i$, $x_i+x_j$ et $x_i-x_j$, 
en  ne gardant que la sous matrice correspondant aux $x_i$),  et les lignes correspondant aux multi-indices $L$ d'ordre 3 (resp 2, resp 1) 
 
  \n {\bf Remarques :}  
  
  $(i)$ On aurait obtenu un r\ésultat analogue en rempla\ç ant les fonctions de 4 variables $x_i.x_j.x_k.x_m$.
 par $x_i+x_j+x_k+x_m$, ou par $(x_i)^2+(x_j)^2+(x_k)^2+(x_m)^2$.
 
  $(ii)$ On n'a pas eu \à utiliser le programme montrant que la courbure \était nulle pour tout $n$, mais seulement pour $n=3$ ou $4$. 
  
 $(iii)$   La m\^eme m\éthode permet de red\émontrer que les tissus ${\cal W}(A_{0,n+3})$ sont tous ordinaires de rang maximum $\pi'\bigl(n,c(n,4)\bigr)$.

  \n {\bf   R\éf\érences   : }
  
 \noindent [B] W. Blaschke et G. Bol, Geometrie der Gewebe, Die Grundlehren der Mathematik 49, Springer, 1938.

  \noindent [Bo] G. Bol, \"Uber ein bemerkenswertes F\"unfgewebe in der Ebene, Abh. Math. Hamburg Univ., 11, 1936, 387-393. 
  
\noindent [CL1] V. Cavalier et D. Lehmann, Ordinary holomorphic webs of codimension one,  arXiv 0703596 v2[mathDS],  2007, et  Ann. Sc. Norm. Super.  Pisa, cl. Sci (5), vol XI (2012), 197-214.



\noindent [GR]  G. Robert, Relations fonctionnelles  polylogarithmiques et tissus plans, pr\épublication  n° 146, 2002,  Universit\é de Bordeaux I.

\noindent [H1] A. H\énaut, On planar web geometry through abelian relations and connections, Ann. of Math. 159 (2004),425-445.


\noindent [Pa]  A. Pantazi, Sur la d\étermination du rang d'un tissu plan, C.R. Acad. Sc. Roumanie 2, 1938, 108-111.

\noindent [Pe]   J. V. Pereira, Resonance webs of  hyperplane arrangements, Advanced studies in Pure Mathematics  99, 2010, 1-30.

\noindent [Pi]  L. Pirio, Sur les tissus planaires de rang maximal et le probl\ème de Chern, note aux CRAS, s\ér. I, 339 (2004), 131-136.

\noindent [PT]  L. Pirio et J.M. Tr\épreau, Abelian functional equations, planar web geometry and polylogarithms, Selecta Mathematica, N.S., 11, n° 3-4, 2005,  453-489. 

 \vskip 2cm

   \n Jean-Paul Dufour, ancien professeur \à l'Universit\é de Montpellier II, \hb
  1 rue du Portalet, 34820 Teyran, France\hb  email : dufourh@netcourrier.com,

  \n Daniel Lehmann, ancien professeur \à l'Universit\é de Montpellier II,\hb   4 rue Becagrun,  30980 Saint Dionisy, France\hb  email : lehm.dan@gmail.com,

    \end{document}